PRELIMINARY VERSION

# POLYGONS IN POLYGONS WITH A TWIST

JAMES M. PARKS

**Abstract**
This is a study of the construction of particular regular sub-*n*-gons *T* in regular *n*-gons *P* using a special system of chords of *P*. The reason for this system is that some of these resulting sub-*n*-gons have areas which are integer divisors *m* of the area of the given *n*-gon *P*. Initially, chords of *n*-gons which determine several well known examples of this relationship were available. However, it has become apparent that a much more general situation exists which encompasses all regular *n*-gons and all integers $m \geq 2$. Dynamic Geometry software was the key to exploring and understanding this relationship.

## 1. INTRODUCTION

Consider the areas of regular polygons and certain of their sub-polygons, which are integer divisors. A particular system of chords is used to determine the sub-polygons. The use of Dynamic Geometry Software [8], [9] is a crucial part of the study.
To help establish the framework for the construction of the system of chords the following 3 examples are presented†.

**Example 1**.
This is a well-known example of a square (a regular *4*-gon), and a particular sub-square (a regular sub-*4*-gon), where the areas are related by the integer multiple $m = 5$, [1], [2], [7, Project 54]. This is a special case of a result called Sylvie's Theorem [3]. There is a simpler case for a triangle [4], but starting with a square is more informative.

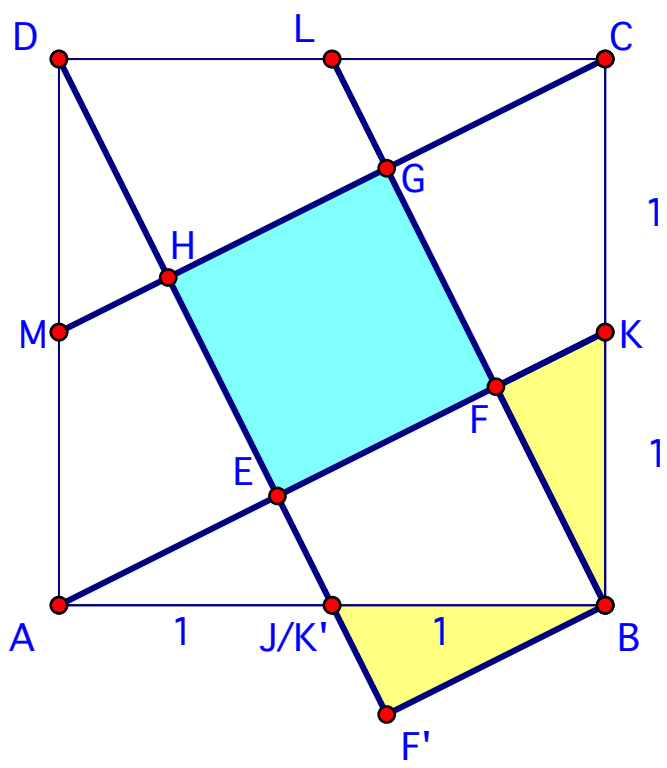

FIGURE 1. Square with Chordal System

Let ABCD be the given square with sides $s = 2$, for convenience, Fig.1. Join each vertex with the midpoint of the opposite side on the right, making a system of 4 chords in a ccw pattern (cw will give equivalent results). The 4 chords intersect to form a sub-square EFGH, by construction (ignore the yellow triangles for now). The claim is the sub-square EFGH has area 1/5 that of the given square ABCD, and thus (ABCD)/(EFGH) = $m = 5$, (parentheses denote area).

†*Sketchpad* was used to construct the examples, and figures in this paper.

This can be argued as follows. Observe that △AJD is a right triangle with legs of length 1 and 2, hence the hypotenuse DJ = √(5). Then, since △AJD and △HMD are similar, and △HMD = △EJA, both results are by construction, we have EJ = 1/√(5), and DH = 2/√(5). So the length of *t*, the side of EFGH, is *t* = EH = DJ - DH - EJ = 2/√(5). Therefore, (EFGH) = $t^2$ = 4/5 units$^2$, and (ABCD)/(EFGH) = 5.
To see this result geometrically, note in the argument above that since BF = FG, and FK = JE = FB/2, the square BFEF', formed by rotating the yellow △BFK ccw 90° about vertex B, is congruent to EFGH. This holds for each of the 4 quadrilaterals on the sides of EFGH, making equal squares on each side of EFGH. Thus, there are 5 equal sub-squares, each of area 4/5 units$^2$, which make up the area of the square ABCD, itself of area 4 units$^2$. Hence (ABCD) = 5(EFGH), so (ABCD)/(EFGH) = 5.

Analytic arguments for the following examples are not given here, but the geometric constructions are.

**Example 2.**
The Hexagon in Fig. 2 has a chordal system with the chords, for example AS, attached from a vertex A to the opposite side ED at a point 1/3 of the distance, ES, on ED. Thus, the perimeter distance *d* from A to S is *d* = 2 1/3 sides. The multiple *m* = 7. This case shows how *Sketchpad* allows one to determine the point S, which determines *d,* by moving point S along ED, until the ratio of the area calculations equal *m*.

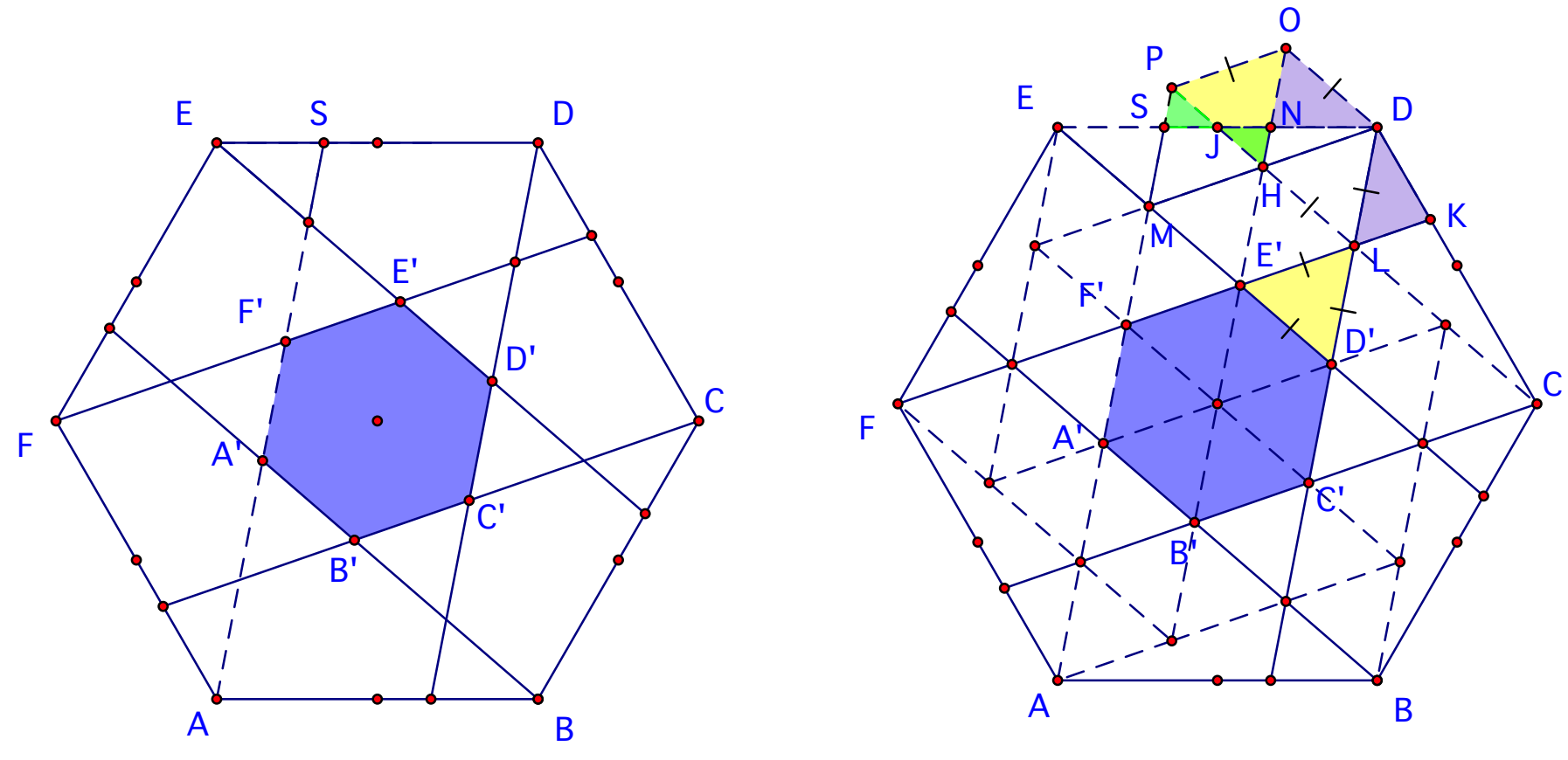

FIGURE 2. Hexagon with Chordal System

For the geometric construction in this case a different argument is used from that given in [6]. Recall that a regular hexagon can be constructed from 6 congruent equilateral triangles.
For this construction, there are 6 “possible” sub-hexagons besides A`…F`, one at each vertex of the sub-hexagon A'…F'. There are 2 types of triangles of interest here, those which are equal to △D'E'L, which is an equilateral triangle, and those which are equal to △KLD. These are in turn associated with each of the “possible” sub-hexagons. If you rotate △KLD about vertex D by an angle of 120°cw, you get △NOD, and it fits on a missing part of the “possible” hexagon centered at H. The other triangle △D'E'L is equivalent to △HPO, by a translation. But the triangle area △JHN is double covered, so rotate one copy of △JHN about J by 180° to coincide with △JPS. This fills in the missing tip of △HPM, and completes the sub-hexagon E'LDOPM. Do this for each of the partial hexagons on each of the vertices of hexagon A'…F'. Thus, there are 6 additional sub-hexagons congruent with the main sub-hexagon A'…F', and this makes 7 sub-hexagons in A…F. So the multiple is *m* = 7.

**Example 3.**
The regular Octagon A…H has the sub-octagon I…P, as shown in Fig. 3, when the chords are identical to AT, where T is the midpoint of the side which is 3 sides away from A in a cw manner, so $d = 2.5$. It is claimed that this chordal system determines the ratio $m = 3$. This can be explained as follows.
The area between the outer octagon A…H, and inner octagon I…P, accounts for 2 more octagons, when split in half as indicated by the 2 black bars in the center figure, Fig. 3. The triangular pieces which make up these two halves are of 3 types. When rearranged as shown in the figure on the right, Fig. 3, they form 2 more octagons. Then the ratio of the hexagon area to the sub-hexagon area is 3, (A…H)/(I…P) = $m$ = 3.

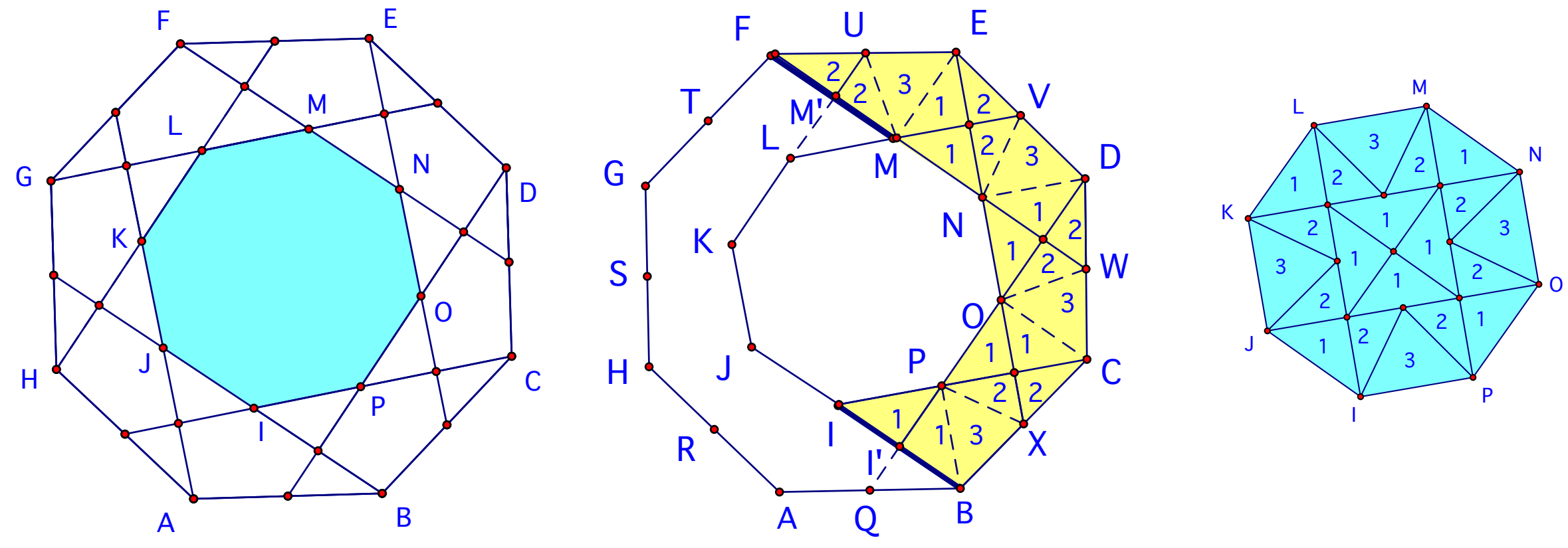

FIGURE 3. Octagon with chordal system

## 2. CONSTRUCTING THE AREA CHORDAL SYSTEM

The notation for the system of chords as used in Example 1. is <1.5>. This means that for a vertex V of the square $S$, a chord $c$ connects V with the midpoint M on the opposite side of $S$, $c$ = VM, which is the distance $d = 1.5$ *sides* away from V, moving along the perimeter, always in a ccw direction, or a cw direction (the results are equivalent). The general definition of the chordal system is given here.

**Definition 1.** *An Area Chordal System for a regular n-gon P, is determined by a chord c, from a vertex A of P to a point S on an opposite side of P from A, c = AS. This chord is then copied to each of the other vertices, by rotation about the center of the n-gon P by the angle 360°/n. The intersections of the chords form a regular sub-n-gon T, by construction. The notation <d> denotes the distance d, in sides, from A to S, which determines the system.*

The examples above may give the impression that the number of examples of this phenomena are few. But if one has a copy of a Dynamic Geometry software available, they will quickly realize that it's the geometric argument which is lacking, while examples abound. This gap will be addressed below.

Regarding the Dynamic Software, since the point S is allowed to moved on a side of $P$, in a chordal system setting, the end points of the other chords will move in unison, the sub-$n$-gon adjusts in size, and any calculations connected to the chord positions, such as area, update automatically. This is Dynamic Geometry in action!

**Definition 2.** *Each Area Chordal System determines an Area Chordal Triple (n, <d>, m). The dimension of the polygon P is n, <d> indicates the chordal system used, and m is the ratio of the area of the n-gon P to that of the sub-n-gon T, (P)/(T) = m.*

The chordal triples for the 3 Examples above are: Square: (4, <1.5>, 5), Hexagon: (6, <2.33…>, 7), and Octagon: (8, <2.5>, 3).
Here are three more examples of known chordal triples: (4, <1.66…>, 13), (6, <2>, 3), and (6, <2.5>, 13).
The chordal triples (6, <2.33…>, 7), and (6, <2.5>, 13) can be found in [3] and [6].

An analytic geometry argument can sometimes be used to compute the area of the sub-$n$-gon $T$, as was done in the Example 1. The *Sketchpad* Polygon Interior Area Calculator is also available to help with some of these calculations.

## 3. AREAS OF REGULAR POLYGONS AND PROPERTIES OF SUB-$n$-GONS

The calculation of the areas of the $n$-gons, using the *Sketchpad* area calculator, is not always as accurate as may be needed or desired. These results can be improved to 5 decimals for all variables by resetting the decimal levels in *Preferences* to their maximums before you begin. This is strongly recommended.

Another way to calculate the area of a regular $n$-gon $P$ is to use the General Area Formula [10].

**Definition 3.** *Let P be a regular n-gon of side s. Then the area of P is* $(P) = ap/2$*, where a is the apothem, thus* $a = s/(2\tan(180°/n))$*, and p is the perimeter,* $p = ns$.

The information available up to this point appears to indicate that only certain regular $n$-gons, and certain chordal systems produce the results on the integer ratio of the areas of the $n$-gons and sub-$n$-gons.
*It will become apparent that this is far from reality.*

A precise value for the side $t$ of the sub-$n$-gon $T$, can be found in terms of the side $s$ and the value of $m$, by the following Proposition.

**Proposition 1.** *If P is a regular n-gon, with side s, T is a regular sub-n-gon with side t determined by a legal chordal system, and* $(P)/(T) = m$*, then* $t = s/\sqrt{m}$.

**Proof.** Consider the equation $(P)/(T) = m$, where $(P) = ap/2$, and $(T) = a'p'/2$, by Definition 3, the General Area Formula for regular polygons. Substituting for $a$ and $a'$, it follows that $(P)/(T) = s^2/t^2$, so $s^2/t^2 = m$, and hence $t = s/\sqrt{m}$.

If $P$ is a unit square, then the calculation of $t$ can also be made in terms of point S on side BC [1, Corollary 4]. The argument is a standard Analytic Geometry exercise, Fig. 1.

**Corollary.** *Let P be a unit square, and T the sub-square determined by the chordal system with point S on side BC. Then* $t = (1-a)/\sqrt{a^2+1}$*, where* $a = |BS|$.

## 4. DETERMINING POINT S ON A SQUARE FOR INTEGER VALUES OF $m$

When working with $P$ a unit square, the calculation of the point S for a particular value of $m$ follows by the Corollary.

**Proposition 2.** *Let P be a unit square with a sub-square T determined by a chordal system for point S on side BC. Suppose T has area an integer divisor m of the area of the P,* $1/(T) = m$*. Then point S is determined as a function of m, by the formula* $|BS| = (m-\sqrt{2m-1})/(m-1)$.

**Proof.** If $m$ is given, then $m = 1/t^2$, by Proposition 1. Since $d = 1+|BS|$, for a point S on BS, $|BS| = (m - \sqrt{2m-1})/(m-1)$, by the Corollary, so point S is determined by $m$.

Finding point S graphically, on an $n$-gon $P$, for a given value of $m$, can at times be a challenge.
For example, if $n = 4$, and $m = 41$, then $|BS| = (41-\sqrt{(81)})/40 = 4/5$, and S = (0.8,1), by the Corollary.
Hence $d = 1.8$.
Another example, let $n = 6$, and $m = 3$, on a unit hexagon $P$. Then $(T) = (3\sqrt{(3)}/2)(1/3) = \sqrt{(3)}/2$, by Proposition 1, using the general area for a hexagon $(3\sqrt{(3)}/2)s^2$ [10]. Also, $(T) = (3\sqrt{(3)}/2)t^2$, so $t = \sqrt{(3)}/3$, $|ES| = 0$, and $S = (0,\sqrt{(3)})$, for $A = (0,0)$, Fig.2. Thus $d = 2$.

The following list is of additional examples, with a few surprises, of chordal triples discovered using these results with the assistance of *Sketchpad*.

| (4, <1.26794…>, 2) | (4, <1.38196…>, 3) | (4, <1.45141…>, 4) | (4, <1.5>, 5) |
|---|---|---|---|
| (4, <1.5366…>, 6) | (4, <1.56574…>, 7) | (4, <1.58957…>, 8) | (4, <1.60961…>, 9) |
| (4, <1.66…>, 13) | (4, <1.75>, 25) | (4, <1.8>, 41) | (4, <1.833…>, 61) |
| (4, <1.875>, 113) | (4, <1.880808…>, 125) | (4, <1.9>, 181) | |
| (5, <1.409…>, 2) | (5, <1.6144…>, 3) | | |
| (6, <2>, 3) | (6, <2.33…>, 7) | (6, <2.5>, 13) | (6, <2.666…>, 31) |
| (6, <2.75>, 57) | | | |
| (7, <1.9328…>, 2) | (7, <2.491…>, 3) | (7, <2.822…>, 4) | |
| (8, <2.1522…>, 2) | (8, <2.5>, 3) | (8, <3.3854…>, 9) | |

The triples (4, <1.75>, 25), (4, <1.88080…>, 125), and (8, <3.3854…>, 9) are discussed in Section 5.

These results lead to the following Theorem for regular polygons and Area Chordal Systems.

**Theorem 1.** *Given a regular n-gon P, and an integer $m \geq 2$, there is a point S on the perimeter of P, such that the Area Chordal System determined by S produces a sub-n-gon T, which satisfies the equation $(P)/(T) = m$.*

**Proof.** Let there be given a regular $n$-gon, say $P$ = A…Z, and an integer $m \geq 2$. Then, by the continuity of the distance function, there is a point S on the perimeter of the $n$-gon $P$, beyond AB (moving ccw), where the chord $c$ = AS, and thus the side distance $d$ from A to S determines the exact integer value of $m$, such that $(P)/(T) = m$. This follows from Proposition 1, since S determines $b$S, $b$ = B or C or …, and thus $T$ in a continuous manner as a distance calculation, and since $(T) = (P)/m$.

## 5. REPLICATION OF EXAMPLES

Consider the case of the replication of a Chordal Triple within itself. Consider the example of the application of the Chordal System (4, <1.5>, 5) to the sub-square EFGH in Example 1, Fig. 4.
This determines a new triple, (4, <1.75>, 25), in the sub-square $T$ = EFGH of the square $S$ = ABCD. Thus there are 25 sub-sub-$n$-gons, name them $U$, each congruent to the square IJKL, and each of area 4/25 units$^2$ in the square ABCD. Therefore, $(S)/(U) = 25$. The triple (4, <1.88080…>, 125) is a further replication of this new triple (4, <1.75>, 25).

The square with the triple (4, <1.5>, 5) is a tiling of the plane, so it's replications are also tilings of the plane [11]. Some polygons are tilings of the plane, but not all Chordal Systems of those which are tiling are also tilings of the plane.

The chordal triple (8, <3.3854…>, 9) is another example of a replication. In this case it is a replication of the triple (8, <2.5>, 3).

Thus, each example of an Area Chordal Triple $(n, <d>, m)$ for a given $n$-gon $P$ can, in theory, be repeated, ad infinitum, in the sub-$n$-gons $T$, with chordal triples $(n, <d_k>, m^k)$, for $k \geq 2$, and $d_k$ the associated distance needed to build the new chordal triple in $P$.

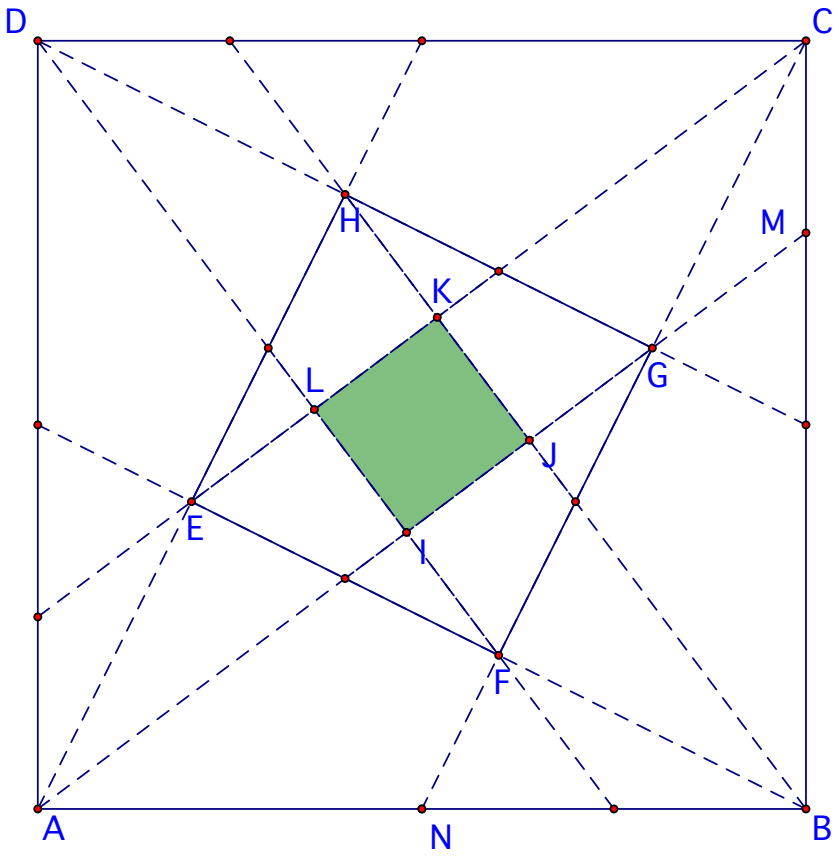

FIGURE 4. Replication of a Square

## 6. SUMMARY

A method of constructing a system of chords of regular polygons $P$ has been discovered which determines regular sub-polygons $T$ of $P$, such that the ratio of the areas of $P$ to $T$ is a positive integer $m \geq 2$, $(P)/(T) = m$. This relationship was originally thought to rely on special choices for chords in order to produce the integer $m$, and was limited to a select few $n$-gons. However, the Dynamic Geometry software demonstrates that there is a technique for finding such sub-polygons $T$ of $P$ such that $(P)/(T) = m$, which applies to *all regular $n$-gons $P$,* and *all integers $m \geq 2$.*

DEPARTMENT OF MATHEMATICS
SUNY POTSDAM
POTSDAM, NY 13676, US
*eMail Address:* parksjm@potsdam.edu